# About some aspects of function interpolation by trigonometric splines

Volodymyr Denysiuk, Dr of Phys-Math. sciences, Professor,

Olena Hryshko

Kiev, Ukraine

National Aviation University

kvomden@nau.edu.ua

**Annotation**

Interpolation of classes of differentiated functions given on a finite interval by trigonometric splines using the phantom node method is considered. This method consists in supplementing a given sequence of values of an approximate function with an even number of values of a phantom function, which is constructed in such a way as to eliminate gaps in both the function itself and its derivatives up to and including a certain order; in the General case, these gaps occur with the periodic continuation of the function given at a finite interval. The results of calculations on test examples for trigonometric splines of the third order are given; these calculations illustrate the high efficiency of the proposed method.

## Introduction

Approximation, respectively, the representation of a known or unknown function $f(x)$ due to the set of some special functions can be considered as a central topic of analysis; such special functions are well defined, easy to calculate, and have certain analytical properties [1]. Algebraic and trigonometric polynomials often act as special functions [2] , polynomial [3,4] and trigonometric [5,6] splines, etc.

An important role in the problems of approximation is played by the accuracy of such an approximation, which is determined by the deviation of the constructed approximating function from a given function $f(x)$. The question of the accuracy of the approximation of both individual functions and their entire classes is studied by the theory of approximations, the foundation of which is laid by the classic works of Chebyshev and Weierstrass, Jackson and Bernstein, and so on.

In the theory of approximations, the problems of approximation content given on classes of functions (or, more generally, on sets of an arbitrary Banach space) are in many cases extremal problems; you need to find the exact upper bound of the approximation error by a given method on a fixed class of functions or specify the best approximation apparatus for this class.

The most significant results of the final nature, ie the results where the solution is brought to exact constants that can no longer be improved, are obtained in the classes of periodic functions, which is easy to explain; classes of periodic functions have a certain symmetry of extreme properties, while the extreme properties of functions given at a finite interval are significantly affected by the disturbing effect of the ends of the interval. Thus, in particular, it is known that on classes of periodic differentiated functions, classes of periodic simple interpolation polynomial splines provide the smallest approximation error in some spaces. [7,8]. The same result can be transferred to classes of simple trigonometric interpolation splines, which coincide with classes of periodic simple interpolation polynomial splines. [9].

In practice, in the vast majority of cases, the classes of differentiated functions are approximated; functions given only at a finite interval have continuous derivatives and are not periodic; therefore, the application of the classical results of the theory of approximations to such functions is impossible. At the same time, the approximation of such functions by polynomial splines is quite simple; more there are estimates of interpolation error that are slightly inferior to estimates of the best approximation [4,10].

The situation is different with simple trigonometric splines; indeed, these splines are intended to approximate only periodic functions and are Fourier series with special coefficients. The approximation of the functions given at a finite interval by Fourier series involves the periodic continuation of these functions over the entire numerical axis; hereinafter the approximate properties of trigonometric splines are determined by the properties of this periodically extended function.

However, with the periodic continuation of the functions given at a finite interval, in the general case there are breaks of the first kind of jump type as the function itself and its derivatives. The presence of such discontinuities sharply degrades the approximate properties of trigonometric Fourier series and, accordingly, trigonometric splines..

Therefore, in our opinion, the task of developing methods for constructing simple trigonometric splines for approximation of classes of differentiated functions given at a finite interval is urgent..

**Formulation of the problem**

Development of methods for approximation by trigonometric splines of classes of differentiated functions given on a finite interval.

**Main part**

In what follows, without losing generality, we will confine ourselves to considering classes of differentiated functions. In the theory of approximations we consider classes of such functions given as on a finite interval $[a,b]$ and $[2\pi]$ - periodic, set on the entire numerical axis. Depending on the conditions satisfied by the senior derivative of these functions, there are different classes of differentiated functions; we will limit ourselves to considering classes $W_{[a,b]}^r$ and $W^r$. Recall that class $W_{[a,b]}^r$ belongs to the set of functions in which the derivative $r-1$ th order locally absolutely continuous on $[a,b]$, a $f^{(r)} \in M[a,b]$, where $M[a,b]$ space measurable, significantly limited to $[a,b]$ functions $f(t)$, norm in which equality is set

$$\|f\|_{M[a,b]} = \sup_{a \le t \le b} \mathrm{vrai} |f(t)|.$$

If $f(t)$ continuous, then $\|f\|_{M(a,b)} = \|f\|_{C[ab]}$. This space, in particular, includes piecewise continuous on $[a,b]$ functions with breaks of the first kind.

In what follows, we will consider these classes of functions on the segment $[0, 2\pi]$.

To denote periodic classes with period $2\pi$, Continuously differentiated functions assigned to all axes, we will still use the same notation without specifying the interval for determining these functions. Therefore, $W^r$ - set set on the entire axis $2\pi$ - periodic functions in which the derivative $r-1$ th order is locally absolutely continuous on any finite interval of the real axis, a $f^{(r)} \in M$, where $M$ space measurable, significantly limited along the entire numerical axis $f(t)$, norm in which equality is set

$$\|f\|_M = \sup_{0 \le t \le 2\pi} \mathrm{vrai} |f(t)|.$$

Let us be given some function $f \in W_{[0,2\pi]}^r$., which we will approximate by trigonometric splines. Since trigonometric splines are trigonometric Fourier series with special coefficients, the problem of presenting this function with trigonometric Fourier series actually arises. It is known that this function $f$ it is necessary to continue periodically on the whole numerical axis. It is clear that with such a periodic continuation in the general case in points $0$ and $2\pi$ there are gaps of the 1st kind of jump type; the presence of such discontinuities significantly impairs the approximate properties of trigonometric Fourier series and, consequently, trigonometric splines. To correct this situation, it is advisable to use the method of phantom functions, which is as follows [11].

**Phantom function method.** By linear replacement of a variable we will display function $f(x) \in W_{M[0,2\pi]}^r$, per segment $t \in [0, 2\pi - \alpha]$, $0 < \alpha < 2\pi$; in between $(2\pi - \alpha, 2\pi)$ define this function so that the periodic continuation of the continuity of the function itself and its derivatives to the order $r-1$, and based on the conditions

$$\lambda(2\pi - \alpha) = f(2\pi); \qquad \lambda(2\pi) = f(0);$$
$$\lambda'(2\pi - \alpha) = f'(2\pi); \qquad \lambda'(2\pi) = f'(0);$$
$$\dots\dots\dots\dots\dots\dots\dots\dots\dots$$
$$\lambda^{(r-1)}(2\pi - \alpha) = f^{(r-1)}(2\pi); \qquad \lambda^{(r-1)}(2\pi) = f^{(r-1)}(0).$$

This is easy to achieve using Hermite interpolation polynomials. Thus, from the consideration of the function $f(x)$, $x \in [a,b]$ we proceeded to consider the function

$$\varphi(t) = \begin{cases} f(t), & t \in [0, 2\pi - \alpha]; \\ \lambda(t), & t \in (2\pi - \alpha, 2\pi). \end{cases}$$

яка, після періодичного продовження з періодом $2\pi$ on the whole numerical axis, is continuous and has continuous derivatives of order $r-1$. It is clear that built in this way periodic with period $2\pi$ is a function of class $W_M^r$, since the added Hermite polynomial does not change the properties of the original function $f(x)$.

Thus, having a function $f(t) \in W_{M[0,2\pi]}^r$ and applying the method of phantom functions, we obtain $2\pi$-periodic function $\varphi(t) \in W_M^r$, which is on the segment $[0, 2\pi - \alpha]$ chalf coincides with the original function.

In many cases, the function under study $f(t) \in W_{M[0,2\pi]}^r$, has a rather complex analytical view or such a view is unknown at all. This function is discretized, ie replaced by a finite sequence of its instantaneous values in the nodes of some, usually equidistant grid, set on $[0, 2\pi]$. In the future there is a problem of interpolation of this sequence by trigonometric splines, which allow to take into account a priori information about the smoothness of the original function $f(t)$. However, in this case it is advisable to use the method of phantom nodes, which is a discrete variant of the method of phantom functions. This method is as follows [11].

**Phantom node method.** Add to the sequence of interpolation nodes on the right side an even number of phantom nodes; values at these phantom points will be set, based, as before, for reasons of crosslinking the function itself and its derivatives. to eliminate the gaps that occur during the periodic elimination of large values of the divided differences around the points $0$ and $2\pi - h$.

It is clear that the addition $2k$ ( $k = 1,2,\dots$ ) phantom nodes increases the number of interpolation nodes on the segment $[0, 2\pi]$, and, accordingly, reduces the step of the interpolation grid, which now becomes equal $h_k = 2\pi(i-1)/(N+2k)$. Since the number of values of the interpolated function does not change, reducing the step of the interpolation grid leads to a decrease in the interpolation segment, which becomes equal to $Nh_k$. На відрізку ж $[2\pi - Nh_k, 2\pi]$ we will build, as before, a function $\lambda(t)$, which satisfies the conditions

$$\lambda(2\pi - Nh_k) = f(2\pi); \qquad \lambda(2\pi) = f(0);$$
$$\lambda'(2\pi - \alpha) = f'(2\pi); \qquad \lambda'(2\pi) = f'(0);$$
$$\dots\dots\dots\dots\dots\dots\dots\dots\dots$$
$$\lambda^{(r-1)}(2\pi - \alpha) = f^{(r-1)}(2\pi); \qquad \lambda^{(r-1)}(2\pi) = f^{(r-1)}(0).$$

As before, this is easy to achieve using Hermite interpolation polynomials. Thus, from the consideration of the function $f(x)$, $x \in [a,b]$ we proceeded to consider the function

$$\varphi(t) = \begin{cases} f(t), & t \in [0, 2\pi - \alpha]; \\ \lambda(t), & t \in (2\pi - \alpha, 2\pi). \end{cases}$$

It is clear that with the periodic continuation of the function $\varphi(t)$ with the period $2\pi$ we get the class function $f(t) \in W_M^r$, for which all the results of the theory of approximations by polynomial splines are applicable [7,8].

Calculating now the value of this function in phantom nodes $t_{N+1}, t_{N+2}, \dots, t_{N+2k-1}$, we obtain the initial data for the construction of trigonometric interpolation splines.

Since the addition of phantom nodes smooths the breaks of the function itself and its derivatives to a certain order, we can expect a reduction in the interpolation error on the segment $[0, 2\pi - Nh_k]$. Note that an even number of points is chosen only so that the total number of interpolation nodes is odd, because interpolation by trigonometric splines of an odd number of points is more convenient. In addition, in our opinion, it is convenient to add a small number of phantom nodes, and put $k = 1, 2$.

In the case where the approximate function is given by its values in the nodes of the uniform grid, when constructing the Hermite polynomial instead of the exact values of the derivatives at points $2\pi - Nh_k$ i $2\pi$ you can use split differences, which are calculated using the step of the new grid applicable applicable $h_k$. Note that the values in phantom nodes can be chosen based on other considerations, thus avoiding the need to build this polynomial.

We illustrate the effectiveness of the proposed method of phantom nodes on the example.

**Example.** Consider the function $f(t) = t + 1$, $t \in [0, 2\pi]$ and put $N = 9$. calculating the value of this function in the nodes of a uniform grid $\Delta_9 = \{t_i\}_{i=1}^{9}$, $t_i = \frac{2\pi}{9}(i-1)$, we get the value $1, 2, 3, 4, 5, 6, 7, 8, 9$. Since in practice simple polynomial cubic splines are used in most cases, we will further limit ourselves to considering simple third-order trigonometric splines. The graph of such a trigonometric interpolation spline that interpolates these values of the function and the graph of the function itself $f(t)$ is shown in Fig.1

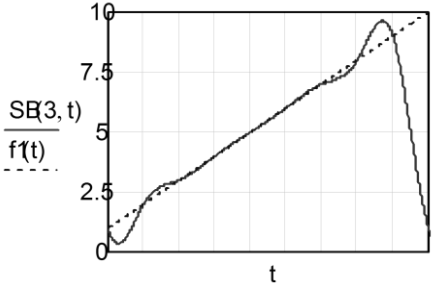

Fig.1. Graphs of interpolation trigonometric spline of order and graph of function $f(t)$.

In the graph shown in Figure 1, the role of the gap that occurs during the periodic continuation of the continuous function is transparent. Although the interpolation error is small in the middle part of the segment, this error is much larger at the ends of the segment.

To the interpolation nodes we will now add phantom nodes, the values of which will be calculated by considering the phantom function $\lambda(t)$ linear; note that the step of the interpolation grid is now $h = 2\pi/11$. A graph of a trigonometric interpolation spline that interpolates these values of the function and a graph of the original function $f(t)$ is shown in Fig.2.

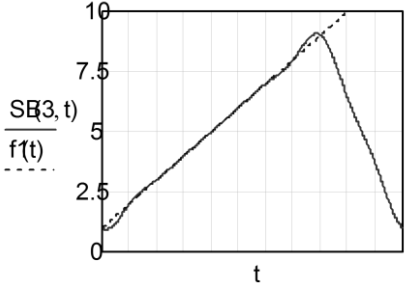

Fig.2. Graph of trigonometric interpolation spline of order and graph of initial function $f(t)$.

The graph shown in Figure 2 illustrates the decrease in error at the ends of the interpolation segment for a linear phantom function.

We now construct a phantom function taking into account the values of the first and second derivatives and calculate the values of this function in two phantom nodes. A graph of a third-order trigonometric interpolation spline interpolating these values of the function and a graph of the original function $f(t)$ is shown in Fig.3.

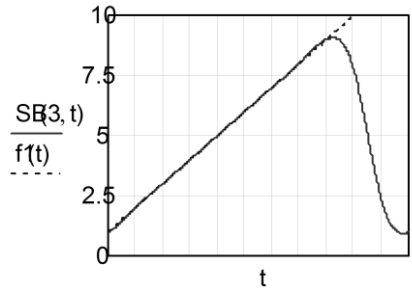

Fig. 3. Graph of trigonometric interpolation spline and graph of initial function $f(t)$.

The graph shown in Figure 3 illustrates the reduction of the error at the ends of the interpolation segment for a phantom function constructed taking into account both the values of the function itself and the values of the first and second derivatives of this function.

The results of further calculations for some other functions are given below.

Table 1.

Function $f(t) = t+1,\ t \in [0, 2\pi)$.

| Number of interpolation grid nodes | Relative interpolation tweak without the use of the phantom node method | Number of added phantom nodes | Error reduction factor Linear phantom function | Error reduction factor Phantom function taking into account the first derivative | Phantom error reduction factor Phantom function with two derivatives |
|---|---|---|---|---|---|
| 9 | .12 | 2 | 3 | 5.5 | 11.6 |
| 9 |  | 4 | 4.3 | 4.3 | 4.3 |

Table 2.

Function $f(t) = t+1,\ t \in [0, 2\pi)$.

| Number of interpolation grid nodes | Relative interpolation tweak without the use of the phantom node method | Number of added phantom nodes | Coefficient reduction factor Linear phantom function | Coefficient reduction factor Phantom function taking into account the first derivative | Phantom tweak reduction factor Phantom function with two derivatives |
|---|---|---|---|---|---|
| 13 | .12 | 2 | 3.2 | 6 | 13.3 |
| 13 |  | 4 | 5 | 5 | 5 |

Table 3.

Function $y = \sin(.75t),\ x \in [0, 2\pi)$.

| Number of interpolation grid nodes | Relative interpolation tweak without the use of the phantom node method | Number of added phantom nodes | Coefficient reduction factor Linear phantom function | Coefficient reduction factor Phantom function taking into account the first derivative | Phantom tweak reduction factor Phantom function with two derivatives |
|---|---|---|---|---|---|
| 9 | .043 | 2 | 4.1 | 31.3 | 45.3 |
| 9 |  | 4 | 2.5 | 9.6 | 40.9 |

Table 4.

Function $y = \sin(.75t),\ x \in [0, 2\pi)$.

| Number of interpolation grid nodes | Relative interpolation tweak without the use of the phantom node method | Number of added phantom nodes | Coefficient reduction factor Linear phantom function | Coefficient reduction factor Phantom function taking into account the first derivative | Phantom tweak reduction factor Phantom function with two derivatives |
|---|---|---|---|---|---|
| 13 | .048 | 2 | 3.8 | 16 | 43.6 |
| 13 |  | 4 | 5.8 | 5 | 171.4 |

Table 5.

Function $y = .02 \exp t,\ x \in [0, 2\pi)$.

| Number of interpolation grid nodes | Relative interpolation tweak | Number of added phantom nodes | Coefficient reduction factor Linear phantom function | Coefficient reduction factor Phantom function taking into account the first derivative | Phantom tweak reduction factor Phantom function with two derivatives |
|---|---|---|---|---|---|
| 9 | .19 | 2 | 1.7 | 4.7 | 2.7 |
| 9 | | 4 | 1.9 | 5.4 | 47 |

Table 6.

Function $y = .02\exp t$, $x \in [0, 2\pi)$.

| Number of interpolation grid nodes | Relative interpolation tweak | Number of added phantom nodes | Tweak reduction factor Лінійна фантомна функція | Phantom error reduction factor c taking into account the first derivative | Phantom error reduction factor c taking into account two derivatives |
|---|---|---|---|---|---|
| 13 | .16 | 2 | 2 | 4.7 | 10 |
| 13 | | 4 | 2.3 | 7.3 | 53.3 |

The results of the calculations shown in Tables 1-6 need to be discussed.

For example, the results of Tables 1-2 show that increasing the number of phantom points from 2 to 4 is impractical; moreover, with such an increase, the number of derivatives taken into account does not reduce the further reduction of the interpolation error. We can assume that this situation arises due to the fact that all derivative functions $f(t) = t + 1$, starting from the second, equal to 0.

For the function $y = \sin(.75t)$ there are exists such values of phantom nodes at which the interpolation error decreases much more than the coefficients shown in tables 3,4. For example, for this function at $N = 9$ there are such values in phantom nodes at which the interpolation error decreases 226.3 times, and at $N = 13$ in the four phantom nodes there are values at which this error is reduced by 1023.8 times. It can be assumed that such anomalously large error reduction coefficients are related to the analyticity of the interpolated function..

For the function $y = .02\exp t$ there are values at which the interpolation error at $N = 13$ decreases by 84.2 times.

It should be noted that the values in phantom nodes at which such interpolation error reduction coefficients are achieved were calculated with an accuracy of .01; by increasing the accuracy of the calculation of these values, it is possible to achieve greater coefficients to reduce the interpolation tweak.

## Conclusions

A method for approximating by trigonometric splines the classes of differentiated functions given on a finite interval is proposed. This method consists in supplementing a given sequence of values of the approximate function with an even number of phantom nodes, the values of which are calculated by sampling the phantom function. In turn, the phantom function is built in such a way as to eliminate breaks in both the function itself and its derivatives to a certain order inclusive; these gaps in the General case occur with the periodic continuation of the functions specified on the finite interval, the entire numerical axis.

The results of calculations on test examples using the proposed methods of trigonometric splines of the third order are given; these calculations illustrate the high efficiency of these methods. Undoubtedly, the proposed method of phantom nodes requires further theoretical research.

## List of references